# Trapezoidal Fuzzy Logic Model for Learning Assessment

**Igor Ya. Subbotin, Professor,**
Department of Mathematics, College of Letters and Sciences,
National University, 5245 Pacific Concourse Drive, Los Angeles, CA 90045, USA
isubboti@nu.edu
Tel. 310 662 2150

**Abstract.** Fuzzy logic has been proven to be extremely feasible in many applications, including engineering, control theory, business, medicine, education, and so on. The current article discusses some applications of fuzzy logic to assessment of learning process. This article is a continuation of the series of the author's papers dedicated to this theme. We consider here a new, more précised fuzzy model for learning assessment.

In 1965, L.A.Zadeh ([Z1], [Z2]) introduced the ideas of so called fuzzy logic as a prospective tool in the control theory for solving some engineering problems that could not be solved with the standard mathematics tools because of their very fuzzy nature. This theory lets us handle and process information in a similar way as the human brain does. Fuzzy logic has been successfully developed by many researchers and has been proven to be extremely productive in many applications (see, for example, [D], [JVR], [KF], [W], [B], and others).

There are some impressive efforts towards the formalizing of the learning process (see, for example, [VS], [VJ], [SMB]). Evaluating our students' work on the on-going basis, we assess their knowledge by assigning grades, labels or pass-non-pass options. Naturally we need to recognize this process has significant fuzziness. While assessing our students' knowledge we are not completely sure about a particular numerical grade or a corresponding label, which could belong to the two adjacent groups of grades with different degrees of membership, which we can measure, for instance, in percentages. Based on commonly used in fuzzy logic ideas, this process could be realized describing so called membership function, which simply assigns to each of the considered element its degree of belonging to the corresponding sets. Formally, it could be described as a function $F = \{(x, f(x): x \in U\}$, where $U$ is the universal set of the discourse, and the range $E(F)$ of function $F$ is $[0, 1]$. It is a very common approach in the USA educational institutions to divide the interval of the specific grades onto three parts and assign the corresponding grade using + and - . For example, $80 – 82 = B-$, $83 – 86 = B$, $87 – 89 = B+$. On the diagram 1 student $Z$ has, let say, 0.25 or 25% degree membership in the set corresponds to the grade C and, therefore, 0.75 or 75% degree membership in the set B, while student Y has the 1.0 degree membership in the set corresponds to the grade B. It follows from the simple geometric considerations that in the described cases degrees of membership for the same element complement each other to 1. It is worthy to note that the same kind of simple geometric arguments imply that this rules are valid for any boundary distributions (for example, 80 - 81 – B - , not 80 – 82 as above, and so on).

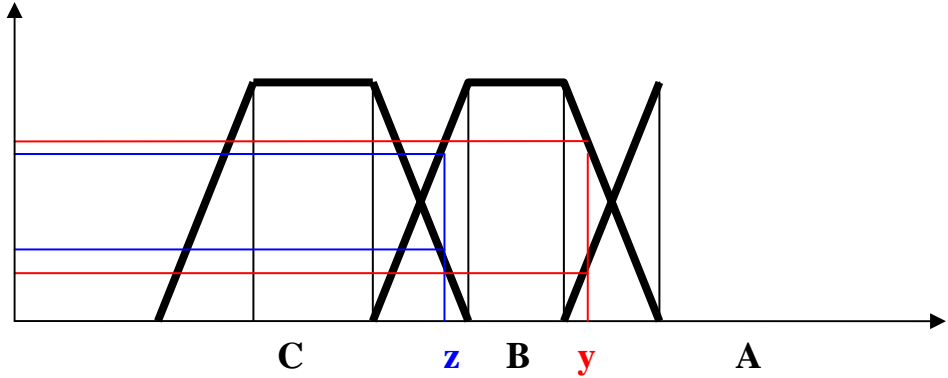

Diagram 1

As some other researchers (see, for example, [VM]), we will base our consideration on the ideas of Voss [VJ], who developed the argument that learning as a specific case of knowledge transfer consists of successive problem-solving activities, in which the input information is represented of existing knowledge with the solution occurring when the input is appropriately represented. This process implements the following states: a) representation of the input data, b) interpretation of this data, c) generalization of the new knowledge, and d) categorization of this knowledge. The states $a$ and $b$ could be unified in one state of interpretation the new knowledge. In the article [VM], the following Fuzzy logic applications have been developed. Let $A_i$, $i = 1,2,3$, be the states of interpretation, generalization, and categorization respectively, and $a,b,c,d,e$ - the linguistic variables of negligible, low, intermediate, high, and complete acquisition of knowledge respectively of each of the $A_i$. Voskoglou considers the set $U = \{a,b,c,d,e\}$ and represents the $A_i$'s as fuzzy sets in $U$. He denotes by $n_{ia}$, $n_{ib}$, $n_{ic}$, $n_{id}$, $n_{ie}$ the numbers of the students that have achieved negligible, low, intermediate, high, and complete acquisition of the state $A_i$ respectively and defines a membership function $m_{A_i}$ by $m_{A_i}(x) = \frac{n_{ix}}{n}$ for each $x \in U$ and, therefore, one can write $A_i = \{(x, \frac{n_{ix}}{n}) : x \in U\}$, where $\sum_{x \in U} m_{A_i}(x) = 1, i = 1,2,3$. A fuzzy relation can be considered here as a fuzzy set of triples, each one of which possess a degree of membership belonging [0, 1]. Consider farther the fuzzy relation

$$R = \{(s, m_R(s) : s = (x, y, z) \in U^3\}$$

where the membership function defined by

$$m_R(s) = m_{A_1}(x)\, m_{R_2}(y) m_{R_3}(z), \text{ for all } s = (x, y, z) \in U^3.$$

This fuzzy relation *R* represents all the possible profiles of student's behavior during the learning process. Further, Voskoglou develops the procedure of comparing few groups of students based on his ideas and supplies the article with examples showing the simplicity of its applications.

We will try to employ another approach to the assessment of students learning. The main base of this approach has been developed in [SB and SBB]. This approach is visible, does not implement any complicated calculations, and, what is important, can be employed to a single student assessment and to the class assessment as well. Depending on evaluation criteria, this approach could be used for the comparing or just for individual independent assessment.

There is a commonly used in the fuzzy logic approach to measure the performance with the pair of numbers $(x_c, y_c)$ as coordinates of the center of mass (the so-called "centroid method") of the represented figure *U*, which we can calculate using the following well-known formulas:

$$(1) \quad x_c = \frac{\iint_F x\,dxdy}{\iint_F dxdy},\ y_c = \frac{\iint_F y\,dxdy}{\iint_F dxdy}.$$

It is not a problem to calculate such numbers using the formulas above; however it could take some significant amount of time. So it would be much more useful in everyday life to simplify the situation as described in diagram 1. For this, we use a rectangular diagram 2. This process implements the following states: C) representation of the input data and interpretation of this data ($y_1$), B) generalization of the new knowledge ($y_2$), and A) categorization of this knowledge ($y_3$). In this case, formulas (1) can be easily transformed to the following simple formulas [SBB]:

$$(2) \quad x_c = \frac{1}{2}\left(\frac{y_1 + 3y_2 + 5y_3}{y_1 + y_2 + y_3}\right),\ y_c = \frac{1}{2}\left(\frac{y_1^2 + y_2^2 + y_3^2}{y_1 + y_2 + y_3}\right).$$

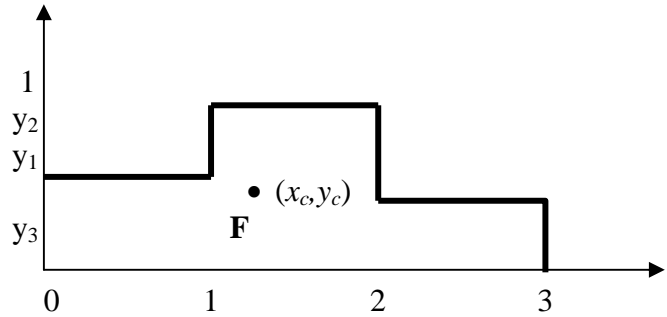

Diagram 2

It is easy to see that the formulas (2) can be generalized for the case when our figure consists not only from three rectangles, but from $n$ rectangles In this case we will come to the following formulas [SBB]:

$$(3) \quad x_c = \frac{1}{2}\left( \frac{\sum_{i=1}^{n}(2i-1)y_i}{\sum_{i=1}^{n} y_i} \right), \quad y_c = \frac{1}{2}\left( \frac{\sum_{i=1}^{n} y_i^2}{\sum_{i=1}^{n} y_i} \right).$$

However, this consideration does not reflect the quite common case when the teacher is not completely sure about the grading and assessing the performance of the students whose performance could be assess as marginal between and close to two adjacent levels. For example, it is something like between 81 and 79 percents. The proposed below "trapezoidal model" helps in this situation.

Instead of rectangles, in the trapezoidal model we use trapezoids. But the most important advantage here is that we allow these triangles have intersections. Following the above mentioned commonly uses grade distribution, we assume that the base of each trapezoid has the length 10, and each adjacent trapezoids have 30% or 3 units of their bases belonging to both of them. This way, we cover the situation of uncertainty of assessment of marginal grades described above.

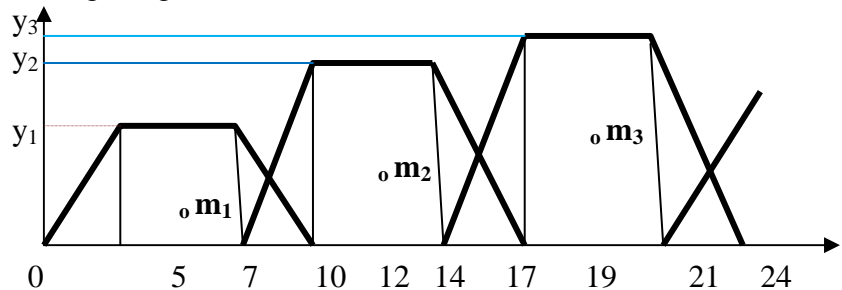

Diagram 3

Note, that the trapezoid form not only better realizes the case of marginal grades. It also correlates the marginal grades through the shape of the triangles in the intersections and this way makes needed correlations.

Since these marginal cases could be considered as common parts for any pair of the adjacent trapezoids, it would be reasonable to represent each trapezoid by its center of mass $\mathbf{m}_i$ and consider the entire figure as the system of these points-centers each of which has its own mass equal to the mass of the corresponding trapezoid.

For the center of mass coordinates we will apply the following commonly used formulas:

$$(4) \quad X_c = \frac{1}{M}\sum_{1}^{n} m_i x_i; \quad Y_c = \frac{1}{M}\sum_{1}^{n} m_i y_i$$

where $m_i$ is the mass of the $i$-trapezoid, and $(x_i, y_i)$ is the coordinates of its center of mass, $M$ is the mass of the considered entire figure. In our specific case, we can assume that $\sum_1^n y_i = 1$. The ordinate Y of the center of mass of a trapezoid could be defined using the formula Y=$\frac{b+2a}{3(a+b)}$ where *b* is the main base, *a* is the upper parallel base. For the sake of simplifying calculations, we assume first that the trapezoids in the Diagram 3 are isosceles, and their bases are 10 units each. Then the abscises of the center of mass of an *i*-trapezoid is  C$x_i$=7*i*-2. Taking into account, that we can allow the density of the figure to be 1, we can equal the mass of a figure to its area. So, $m_i = S_i = 7y_i$. Taking into account that $\sum_1^n y_i = 1$, after some elementary geometric reasoning we find that the coordinates of the center of mass of the *i*-trapezoid are

$$Cx_i = 7i\text{-}2$$
$$Cy_i = \frac{3}{7}y_i,$$

while its area $S_i = 7y_i$.

Therefore by (4) we obtain

$$M = 7\sum_1^n y_i = 7.$$

(5) $\quad X_c = \frac{1}{M}\sum_1^n 7y_i(7i-2) = \frac{\sum_1^n 7y_i(7i-2)}{7\sum_1^n y_i} = \frac{\sum_1^n y_i(7i-2)}{\sum_1^n y_i} = 7\sum_1^n y_i i - 2;$

$$Y_c = \frac{1}{M}\sum_1^n 7y_i \frac{3}{7}y_i = \frac{3\sum_1^n y_i^2}{7\sum_1^n y_i} = \frac{3\sum_1^n y_i^2}{7}.$$

In order to use the formulas in the case when the bigger base b of a trapezoid is 1, we have

$$M = 0.7\sum_1^n y_i = 0.7.$$

(6) $\quad X_c = \frac{1}{M}\sum_1^n 0.7y_i(7i-2) = \frac{\sum_1^n 0.7y_i(0.7i-0.2)}{0.7\sum_1^n y_i} = \frac{\sum_1^n y_i(7i-2)}{\sum_1^n y_i} = 0.7\sum_1^n y_i i - 0.2;$

$$Y_c = \frac{1}{M}\sum_1^n 0.7y_i \frac{3}{7}y_i = \frac{3\sum_1^n y_i^2}{7\sum_1^n y_i} = \frac{3\sum_1^n y_i^2}{7}.$$

Now consider some example.

In the article [SBB] we have considered some examples of comparing the performances of two classes in some marginal cases. In the case when it is not clear how to decide which class performance is better we usually compare classes' GPAs (Grade Point Averages in the American system's meaning), and the indicators commonly called "the quality of knowledge" – the ratio of the sum of the numbers of all B and A to whole amount of grades. There are some ambiguous cases here. For example, consider the following two classes' grades:

| Ratio of the class students reached the following stage of knowledge acquisition | Class I | Class II |
| --- | --- | --- |
| C | 10 | 0 |
| B | 0 | 20 |
| A | 50 | 40 |

For these both classes the GPA is 3.7. "The quality of knowledge" for the second class is higher than for the first one. The standard deviation for the second class is definitely smaller. So from the common point of view and from the statistical point of view the situation in the second class is better. However, some instructors could prefer the situation in the first class, since there are much more "perfect" students in this class. Everything is determined by the set of goal preferences. Using formulas (3) for the case $n=3$ we obtaine.

$$x_{c1} = \frac{1}{2}\left(\frac{1}{6} + 3 \cdot \frac{0}{6} + 5 \cdot \frac{5}{6}\right) = \frac{13}{6} = 2\frac{1}{3}, \quad y_{c1} = \frac{1}{2}\left(\frac{1}{36} + 0^2 + \frac{25}{36}\right) = \frac{13}{36}.$$

and

$$x_{c2} = \frac{1}{2}\left(0 + 3 \cdot \frac{1}{3} + 5 \cdot \frac{2}{3}\right) = \frac{13}{6} = 2\frac{1}{3}, \quad y_{c2} = \frac{1}{2}\left(0^2 + \frac{1}{9} + \frac{4}{9}\right) = \frac{5}{18}.$$

As we can see, $x_{c1} = x_{c2} = 2\frac{1}{3} = 2.33$ and $y_{c1} > y_{c2}$. It means that in this case, considering the described above rectangular model, we find that the centers of mass lie on the same vertical line $x = 2.33 > 1.5$, and the first center is a bit higher. So the first class performs better by the following standards from [SBB].

*Among two or more classes the class with the biggest $x_c$ performs better; If two or more classes have the same $x_c \geq 1.5$, then the class with the higher $y_c$ performs better. If two or more classes have the same $x_c \leq 1.5$, then the class with the lower $y_c$ performs better.*

Consider the same example using our formulas (6) for the trapezoidal model.

| Y | Class I | Class II |
|---|---|---|
| $Y_1$ | 1/6 | 0 |
| $Y_2$ | 0 | 1/3 |
| $Y_3$ | 5/6 | 2/3 |

We obtain:
$$X_{c1} = 5/3, \; Y_{c1} = 13/42,$$
$$X_{c2} = 5/3, \; Y_{c2} = 10/42$$

So, our trapezoidal model also confirms that the second class performs slightly better than the first one.